\documentstyle[aps,epsfig,amssymb,balanced,times,prl]{revtex}
\begin{document}
\title{Tails of probability density for sums of random independent variables}
\author{Michael I. Tribelsky\cite{E-mail}}
\address{Department of Applied Physics, Faculty of Engineering, Fukui
University, Bunkyo 3-9-1, Fukui 910-8507, Japan}
\date{\today}
\maketitle
\begin{abstract}
The {\it exact} expression for the probability density $p_{_N}(x)$
for sums of a {\it finite} number $N$ of random independent terms
is obtained. It is shown that the very tail of $p_{_N}(x)$ has a
Gaussian form if and only if all the random terms are distributed
according to the Gauss Law. In all other cases the tail for
$p_{_N}(x)$ differs from the Gaussian. If the variances of random
terms diverge the non-Gaussian tail is related to a L\'{e}vy
distribution for $p_{_N}(x)$. However, the tail is not Gaussian
even if the variances are finite. In the latter case $p_{_N}(x)$
has two different asymptotics. At small and moderate values of $x$
the distribution is Gaussian. At large $x$ the non-Gaussian tail
arises. The crossover between the two asymptotics occurs at $x$
proportional to $N$. For this reason the non-Gaussian tail exists
at {\it finite\/} $N$ only. In the limit $N$ tends to infinity the
origin of the tail is shifted to infinity, i. e., the tail
vanishes. Depending on the particular type of the distribution of
the random terms the non-Gaussian tail may decay either slower
than the Gaussian, or faster than it. A number of particular
examples is discussed in detail.
\end{abstract}
\pacs{PACS numbers:  02.50.Ey, 
05.40.-a, 
05.40.Fb, 
05.45.Tp 
}
\begin{twocolumns}
It is well known that if a random quantity $x$ is a sum of a
number of independent random variables $\xi$
\begin{equation}\label{sum}
 x=\sum_{n=1}^{N}\xi_n,
\end{equation}
the corresponding probability density $p(x)$ is described by a
Gaussian distribution (GD), provided the conditions of the Central
Limit Theorem (CLT) of Probability Theory\cite{GK} are fulfilled.
However, in reality while GD fits well the profile of $p(x)$ at
small and moderate vales of $x$ the very tail of $p(x)$ often
deviates from the Gauss Law. Most of the observed deviations are
related to the so-called {\it heavy tails}\cite{SF,war,Stan,turb},
which exhibit a decay slower than the Gaussian (usually the heavy
tails decay as a certain power of $x$). It should be stressed that
despite the heavy tails do not influence much the average
characteristics of $x$, such as, e.g., differen moments of it,
large fluctuations of $x$ are described by the tails entirely.
Therefore the tails of $p(x)$ play an important role in those
problems, where rare but large fluctuations are of special
interest, for example for temperature fluctuations in a nuclear
reactor, where too big fluctuations may result in damage of the
reactor.

The heavy tails have been observed in a wide diversity of
phenomena ranging from the statistics of solar flares\cite{SF} to
the one of casualties of wars\cite{war} and from market price
fluctuations\cite{Stan} to fluid dynamics\cite{turb}. It provides
grounds to suppose that the phenomenon should have a general
cause, which is not related to a particular realization of a
random sum in a given problem. However, to the best of my
knowledge such a cause has not been revealed yet.

The present contribution is an attempt to fill the gap. It is
shown that occurrence of non-Gaussian tails is quite a common
case, which does not require any special conditions to be
realized. Visa versa, special (and very strict) conditions are
required for realization of Gaussian tails. It is shown also that
depending on the distribution of random $\xi$ the tails of $p(x)$
could be either heavy, or {\it light\/}. In the latter case the
tails decay faster than the Gaussian. Examples of tails of both
types are considered in detail.

For the sake of simplicity in what follows sums of independently
and identically distributed random variables with zero mean values
are considered. Extension of the approach to more general cases is
straightforward.

First of all notice that, of course, nothing is wrong with CLT.
The point is that CLT says that $p(x)$ converges to GD {\it
asymptotically\/} at $N \rightarrow \infty$. On the other hand in
any practical case $N$ is a {\it finite\/} quantity. At finite $N$
the profile of probability density for $x$ may have nontrivial
dependence on $N$. To stress this dependence $p(x)$ at finite $N$
will be denoted below as $p_{_N}(x)$. Then, if $p_{_N}(x)$ has a
non-Gaussian tail, which begins at $|x| \gg x_c(N)$ and if $x_c(N)
\rightarrow \infty$ at $N \rightarrow \infty$, it does not
contradict CLT. What I am going to show is that the occurrence of
such a tail is a generic feature of the problem.

To begin with let us consider two following examples
\begin{equation}\label{step}
  f(\xi) = \left\{ \begin{array}{l} 0 \; {\rm at} \; \xi < -1/2\\
  1\; {\rm at}\; -1/2 \leq \xi \leq 1/2\\
  0\; {\rm at} \;\xi > 1/2\end{array} \right.
\end{equation}
 and
\begin{equation}\label{m}
  f(\xi) \cong A/|\xi|^m \; {\rm at}\; |\xi| \gg \xi_c,
\end{equation}
where $f(\xi)$ denotes the probability density for random $\xi$.
In Eq.~(\ref{m}) $A,\; m$ and $\xi_c$ stand for certain positive
constants.

In case of Eq.~(\ref{step}) the conditions of CLT hold. However,
it is evident that at any finite $N\;p_{_N}(x)$ vanishes
identically at $|x| \geq N/2$ instead of exhibiting the
exponential Gaussian decay at $x \rightarrow \infty$\cite{Tak}. In
this case one has a {\it superlight\/} tail (actually, no tail at
all).

To see a non-Gaussian tail in case of Eq.~(\ref{m}) let us first
calculate the probability $p_1(x)dx$ that $x$ has a certain large
value from the interval $x_1 \leq x \leq x_1+dx$ because of
contribution of a single, anomalously large term, while all other
terms in sum Eq.~(\ref{sum}) are much smaller than this single
term. The probability that a term in Eq.~(\ref{sum}) has the value
$\xi_n=x$, where $x$ belongs to the specified region is $f(x)dx$.
The probability, that $|\xi_n|$ is much smaller then $|x|$ is
$\int_{-\xi}^{\xi}f(\xi^\prime)d\xi^\prime,\; \xi \ll |x|$. Due to
independence of $\xi_n$ the probability that, say, $\xi_1=x$ and
$|\xi_n| \ll |x|$ at $2 \leq n \leq N$ is
$[\int_{-\xi}^{\xi}f(\xi^\prime)d\xi^\prime]^{N-1}f(x)dx$.
Finally, bearing in mind that there are $N$ terms in sum
Eq.~(\ref{sum}) and every one may be large, we obtain the
following expression for the total probability of the desired
event
\begin{eqnarray*}
 & & p_1(x)dx=
  \left[\int_{-\xi}^{\xi}f(\xi^\prime)d\xi^\prime\right]^{N-1}Nf(x)dx\equiv\\
  & &\left[1- \int_{-\infty}^{-\xi}f(\xi^\prime)d\xi^\prime -
  \int_{\xi}^{\infty}f(\xi^\prime)d\xi^\prime\right]^{N-1}Nf(x)dx.
\end{eqnarray*}
Let us suppose that $x$ is so large that the following conditions
hold
\[  N\!\left(\int_{\xi}^{\infty}\!f(\xi^{\prime})d\xi^\prime
+  \int_{-\infty}^{-\xi}\!f(\xi^{\prime})d\xi^\prime\right)\! \ll
1; \;\;
  \xi_c \ll \xi \ll x. \]
In this case the probability density $p_1(x)$ is reduced to the
expression
\begin{equation}\label{p1as}
  p_1(x) \cong Nf(x),
\end{equation}
i.e., a heavy tail arises, see Eq. (\ref{m}).

Note, that according to CLT the probability density $p(z)$ for the
normalized sum
\[z\equiv x/(\sigma \sqrt N),\]
where $\sigma^2\equiv \langle \xi^2 \rangle$,
should converge to GD at $m>3$, viz.
\begin{equation}\label{GD}
  p(z) \cong p_{_G}(z) \equiv \frac{1}{\sqrt{2\pi}}
  \exp\left(-\frac{z^2}{2}\right)\; {\rm at}\; N \rightarrow
  \infty.
\end{equation}
Here and in what follows $\langle \ldots \rangle$ denotes average
over $\xi$.

GD Eq.~(\ref{GD}) describes the probability density for the sum to
have a value $x$ due to contribution of a big number of individual
terms. Naturally, it cannot be smaller than the probability
density related to contribution of a single anomalously big term.
In other words GD Eq.~(\ref{GD}) certainly fails to describe the
actual situation, when the condition $p_{_G}(z)\geq
p_1(x(z))dx/dz$ is violated\cite{note}. It provides us with the
applicability condition for the GD $|z| \ll z_{_G}$, where
$z_{_G}$ stands for the solution of the equation $p_{_G}(z)=
p_1(x(z))dx/dz$. At large $N$ the approximate solution of this
equation obtained by iterations is in the following
\begin{eqnarray*}
  z_{_G} &\approx &\left[(m-3)\ln N + \right. \\ 
  & & \left. m \ln\left((m-3)\ln N +
  \ldots\right)\right]^{1/2} \gg 1.
\end{eqnarray*}
Thus, for the problem in question the GD is invalid when the value
of $x$ becomes much greater than the standard deviation for the
given GD.

In what follows I am going to prove that the above result is a
generic property of the problem, which is not connected with the
particular type of $f(\xi)$. Naturally, it is supposed that
$f(\xi)$ decays at big $|\xi|$ sharp enough, so that the variance
$\sigma^2$ remains finite. Otherwise instead of GD one will obtain
a L\'{e}vy distribution (LD). The case of LD will be considered
later separately.

Let us calculate the characteristic function $\varphi_{_N}(\omega)
\equiv \langle \exp(i\omega z)\rangle$, where $\omega$ is a real
quantity. By definition
\begin{equation}\label{char}
  \varphi_{_N}(\omega) \equiv
  \int_{-\infty}^{\infty}\exp(i\omega z)p_{_N}(z)dz
\end{equation}
On the other hand, employing Eq. (\ref{sum}) and independence of
$\xi_n$ one can obtain
\begin{eqnarray*}
& &\varphi_{_N}(\omega) \equiv \\
& & \left\langle {\exp \left[\sum\limits_{n = 1}^N
{\frac{{i\omega\xi _n }} {{\sigma \sqrt N }} } \right]}
\right\rangle \equiv \left\langle {\prod\limits_{n = 1}^N {\exp
\left[\frac{{i\omega \xi _n}} {{\sigma \sqrt N }} \right]}}
\right\rangle \equiv
\\
& & \prod\limits_{n = 1}^N {\left\langle {\exp
\left[\frac{{i\omega \xi _n}}{{\sigma \sqrt N }} \right]}
\right\rangle } \equiv \left\langle {\exp \left[\frac{{i\omega
\xi}} {{\sigma \sqrt N }} \right]} \right\rangle ^N  \equiv
\left[{g_{_N} (\omega )} \right]^N,
\end{eqnarray*}
where
\begin{equation}\label{gN}
  g_{_N}(\omega) \equiv \int\limits_{ - \infty }^\infty
  {\exp \left[\frac{{i\omega \xi}}{{\sigma \sqrt N }}\right]} f(\xi
  )d\xi.
\end{equation}
Eqs. (\ref{char}) -- (\ref{gN}) supplemented with the inverse
Fourier transform give rise to the following expression for
$p_{_N}(z)$
\begin{equation}\label{p}
  p_{_N}(z ) = \frac{1} {{2\pi }}\int\limits_{ - \infty }^\infty
{\exp [ - i\omega z  + N\ln g_{_N} (\omega )]d\omega }.
\end{equation}
It should be stressed that Eqs. (\ref{char}) -- (\ref{p}) are just
a sequence of {\it identities}. Thus, expression Eq.~(\ref{p}) is
an {\it exact\/} result, valid at {\it any} value of $N$\cite{n}
\begin{figure}
\vspace*{-8cm}
\epsfxsize = 0.65\textwidth 
\makebox{\hspace{-1.8cm} \epsfbox{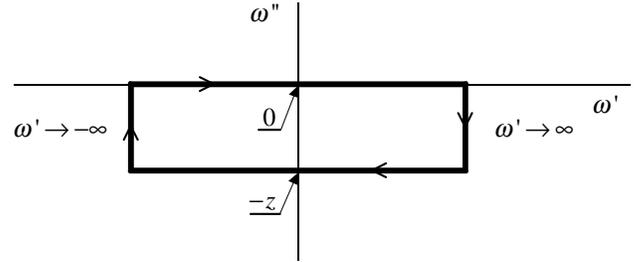}} \vspace*{-4cm}
\caption{Contour of integration for \protect{Eq.~(\ref{pG})}
\label{F1}}
\end{figure}

To obtain GD from Eqs.~(\ref{gN}) -- (\ref{p}) let us consider the
limit of large $N$ and small $\omega$. Expansion of the integrand
in Eq.~(\ref{gN}) in powers of $\omega/\sqrt N$ and integration (I
remind the reader that $\langle \xi \rangle = 0$) yield
\begin{eqnarray}
  & & g_{_N}(\omega) = 1 - \frac{\omega^2}{2N} +
  o\left(\frac{\omega^2}{2N}\right)\label{gG}\\
  & & N\ln g_{_N}(\omega)=-\frac{\omega^2}{2} + o(\omega^2)\label{lngN}
\end{eqnarray}
Then, inserting Eq.~(\ref{lngN}) in Eq.~(\ref{p}) one obtains
\begin{equation}\label{pG}
  p_{_N}(z ) \cong \frac{1} {2\pi }\int\limits_{ - \infty }^\infty \exp \left[
- \frac{1} {2}(\omega  + iz )^2  - \frac{z ^2 } {2} \right]
d\omega.
\end{equation}
The integral is taken by integration in the plane of complex
$\omega$ along the contour shown in Fig.~\ref{F1}.

The integrand is an analytic function inside the contour, hence
the integral over the entire contour is zero. The integrals along
the vertical segments vanish exponentially at $\omega^\prime
\rightarrow \pm \infty$. Therefore the integral Eq.~(\ref{pG})
equals the one along the line $-\infty - iz < \omega < \infty -
iz$. It eventually results in GD Eq.~(\ref{GD}).

However, here we are focusing on the applicability conditions of
these transformations rather then on the well-known result of the
integration. The transformation of the integral Eq.~(\ref{pG})
into the one in the complex plane implies that expansion
Eq.~(\ref{gG}) is now performed in the complex plane too. Thus,
the applicability conditions for the expansion say that
$|\omega|^2/N \ll 1$, where $\omega \equiv \omega^{\prime} +
i\omega^{\prime \prime}$. This condition must hold at least for
those $\omega$'s, which make the main contribution to the
integral. For the integral along the line $-\infty - iz < \omega <
\infty - iz$ the main contribution is made by $\omega^\prime =
O(1)$. However, the imaginary part of these $\omega$'s is $-z$. We
are interested in the tail of $p_{_N}(z)$, where $|z| \gg 1$, so
for $\omega$ with $\omega^\prime = O(1) \ll |z|$ we have $|\omega|
\cong |z|$, which finally yields the following applicability
condition for GD Eq.~(\ref{GD})
\begin{equation}\label{GDappl}
  z^2 \ll N
\end{equation}

Generally speaking, it means that the tail of actual $p_{_N}(z)$
given by Eq.~(\ref{p}) differs from that for Eq.~(\ref{GD}) at
$z^2$ of order of (or grater than) $N$. The only exception from
the rule is the case when
\begin{equation}\label{gGeq}
  g_{_N}(\omega) \equiv \exp \left(-\frac{\omega^2}{2N}\right),
\end{equation}
and expansion Eq.~(\ref{lngN}) transforms into an identity.
Equalizing Eq.~(\ref{gN}) and Eq.~(\ref{gGeq}) and employing the
inverse Fourier transform, one can recover the only $f(\xi)$,
which brings about a Gaussian tail for $p_{_N}(x)$. This $f(\xi)$
is as follows
\[
f(\xi)=\frac{1}{\sigma\sqrt{2\pi}}
\exp\left(-\frac{\xi^2}{2\sigma^2}\right).
\]
It is worth mentioning that contrary to all other cases, when GD
is an asymptotic law valid at large $N$ only, now Eq.~(\ref{GD})
is the exact result, which is good at any $N$ and any $z$. This
property reflects the well known fact of self-similarity of a
Gauss distribution\cite{Zas}

Thus, we have arrived at a number of important conclusions
\begin{itemize}
  \item A sum of {\it any finite\/} number of independently and
  identically distributed random variables has a Gaussian tail
  {\it if and only if\/} the random variables themselves have a
  Gaussian distribution.
  \item In any other situation the probability density
  $p_{_N}(z)$ for the normalized sum $z$ has a tail, which differs
  from the Gaussian.
  \item Contrary to the GD for the normalized sum $z$ the
  non-Gaussian tail is not universal. Its profile depends on the
  distribution of the random terms contributing to the sum and may
  be either heavier than the Gaussian, or lighter than it.
  \item The non-Gaussian tail begins at $|z|$ of order of $\sqrt
  N$ and lasts up to infinity.
  \item The non-Gaussian tails exist at finite $N$ only. It
  disappears in the limit $N \rightarrow \infty$.
\end{itemize}

To illustrate these general conclusions let us consider in detail
several particular examples of $f(\xi)$ and obtain the
corresponding profiles of $p_{_N}(z)$ based upon Eqs.~(\ref{gN})
-- (\ref{p}) . In case of Eq.~(\ref{step}) simple calculations
yield
\[ g_{_N}(\omega) = \frac{2\sigma\sqrt N}{\omega}\sin\frac{\omega}
{2\sigma\sqrt N}\:, \;\; \sigma^2 = \frac{1}{12},\]
so that
\begin{equation}\label{pstep}
  p_{_N}(z) = \frac{1}{2\pi}\int_{-\infty}^{\infty}
  \left(\frac{2\sigma\sqrt N}{\omega}\sin\frac{\omega}
  {2\sigma\sqrt N}\right)^N e^{-i\omega z}d\omega.
\end{equation}

First I show that the above expression for $p_{_N}(z)$ does
fulfill the condition $p_{_N}(z) \equiv 0$ at $|z| \geq \sqrt
N/(2\sigma)$. To this end let us extend the integration in
Eq.~(\ref{pstep}) to the plane of complex $\omega$ and note that
the integrand is an analytic function of $\omega$ at any finite
$|\omega|$. For this reason the corresponding integral over any
closed contour in the complex plane is identical zero. For
definiteness in what follows I suppose that $z>0$ (the case of
negative $z$ is analyzed analogously). Then, I consider a closed
contour consisting of two elements --- a segment of the real axis
$-R \leq \omega^\prime \leq R$ and a circular arc with the radius
$R$ lying in the lower half-plane, which connects the edges of the
segment. Next, note that the greatest term of the integrand in the
lower half-plane has the form
\[\frac{{\rm const}}{\omega^N}\exp\left[i\omega\left(
\frac{\sqrt N}{2\sigma}-z\right)\right] \] It is straightforwardly
seen from the above expression that the integral over the arc
tends to zero at $z > \sqrt N/(2\sigma)$ and $R \rightarrow
\infty$. Finally, bearing in mind that the integral over the
entire closed contour is zero, I conclude that the integral
Eq.~(\ref{pstep}) does turns into identical zero at $z > \sqrt
N/(2\sigma)$. As an example the plot $p_{_{N}}(z)$, obtained by
numerical integration of Eq.~(\ref{pstep}) at $N=10$ is presented
in Fig.~(\ref{F2}).

Next, we discuss $p_{_N}(z)$ at $f(\xi)$ resulting in heavy tails.
Specifically we consider the normalized probability density
$f(\xi)$ in the form
\begin{equation}\label{fl}
  f(\xi) = \frac{A}{1 + \xi^{2l}}, \;\;
  A \equiv \frac{l}{\pi}\sin \frac{\pi}{2l},
\end{equation}
where $l$ is any positive integer number. Thus, the case in
question is a particular realization of $f(\xi)$ given by
Eq.~(\ref{m}). Straightforward calculation yields the following
expression for the variance at $l \geq 2$
\[  \sigma^2 = \left(\sin\frac{\pi}{2l}\right)/
\sin\frac{3\pi}{2l}.
\]

At $l=1$ the variance diverges. It is trivial to obtain $p(z)$ at
$l=1$. The only modification to the described approach is that the
normalized sum $z$ in this case should be defined as $x/N$. For
this $z$ Eq.~(\ref{gN}) yields $g_{_N}(\omega)=\exp(-|\omega|/N)$.
Next, Eq.~(\ref{p}) brings about the final formula (a L\'{e}vy
distribution)
\begin{equation}\label{Levy}
  p(z)= \frac{1}{\pi} \frac{1}{1 + z^2}.
\end{equation}
It should be stressed (i) the expression Eq.~(\ref{Levy}) is an
exact result, valid at any $N$, and (ii) the profile of $p(z)$ in
Eq.~(\ref{Levy}) coincides with that for $f(\xi)$ in
Eq.~(\ref{fl}) at $l=1$ --- a L\'{e}vy distribution possesses the
same self-similar properties as a Gaussian does\cite{Zas}.
However, here we are focused on the tail behavior rather then on
the self-similarity of distributions. Following this line, I would
like to emphasize that the tail of $p_{_N}(x)$, which
Eq.~(\ref{Levy}) leads to, has the form of Eq.~(\ref{p1as}) and
obviously may be explained in the same terms\cite{note}.
\begin{figure}
\epsfxsize = 0.45\textwidth 
\makebox{
\epsfbox{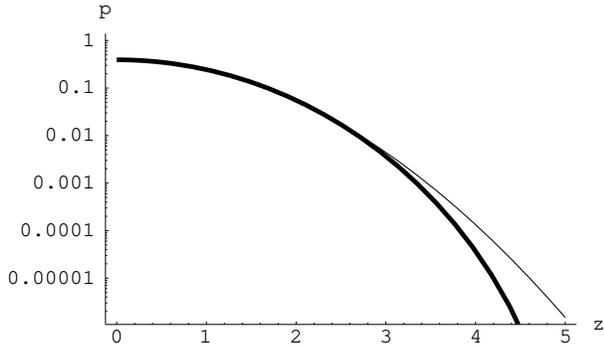}} \vspace*{5mm} \caption{\protect{Plots of
Gaussian distribution Eq.~(\ref{GD}) for the normalized sum $z$
(thin line) and probability density $p_{_{N}}(z)$ calculated
numerically according to Eq.~(\ref{pstep}) at $N=10$ (thick line)
exhibiting a tail {\it lighter\/} than the Gaussian. Note, that in
agreement with Eq.~(\ref{GDappl}) deviation of $p_{_{N}}(z)$ from
the GD begins at $z \approx \sqrt N$. }\label{F2}}
\end{figure}

To calculate $p_{_N}(z)$ in case of arbitrary $l \leq 2$ let us
employ the fact that $g_{_N}(\omega)$ is a real even function of
$\omega$ at any real even $f(\xi)$, see Eq.~(\ref{gN}), and hence
for such $f(\xi)$ Eq.~(\ref{p}) may be transformed as follows
\begin{equation}\label{pl>2}
  p_{_N}(z) = \frac{1}{\pi}{\rm Re}\int_0^\infty
  \left(g_{_N}(\omega)\right)^N e^{-i\omega z}dz.
\end{equation}
Thus, it suffices to calculate $g_{_N}(\omega)$ for positive
$\omega$ only. Next, extending integral Eq.~(\ref{gN}) into the
plane of complex $\xi$, note that the for $f(\xi)$ given by
Eq.~(\ref{fl}) integrand has $l$ simple poles in the upper
half-plane. The poles are located in the points
\[ \xi_j = \exp \left[ \frac{i\pi}{2l}(2j+1)\right], \;\; 0 \leq j \leq l-1. \]
Then, the integral is equal the following expression
\begin{equation}\label{gsum}
g_{_N}(\omega) =
    i\left( \sin \frac{\pi}{2l} \right)\sum_{j=o}^{l-1}
    \frac{\exp\left[ \frac{i\omega}{\sigma \sqrt N}
    \exp\left[\frac{i\pi}{2l}(2j+1)\right]\right]}
    {\exp\left[\frac{i\pi}{2l}(2j+1)(2l-1)\right] }
\end{equation}
Calculation of $p(z)$ according to Eqs.~(\ref{pl>2}), (\ref{gsum})
is rather cumbersome, see Appendix. It yields the following result

\begin{equation}\label{plf}
  p_{_N}(z) \cong \frac{lN}{\pi(\sigma \sqrt
  N)^{2l-1}}\frac{\sin \frac{\pi}{2l}}{z^{2l}}
  \;\;\; \mbox{at}\;\;\; z \gg \frac{\sqrt N}{\sigma}\:.
\end{equation}

Plots of $p_{_N}(z)$ obtained by numerical integration of
Eqs.~(\ref{pl>2}), (\ref{gsum}) at $l=2$ and two values of $N$
($N=10^2,\; 10^4$, respectively) are presented in Fig.~{\ref{F3}}.
The crossover from GD to heavy tails Eq.~(\ref{plf}) is seen
clearly. In agreement with the above discussion the larger is $N$
the longer the corresponding curve follows GD.
\begin{figure}
\epsfxsize = 0.45\textwidth 
\makebox{
\epsfbox{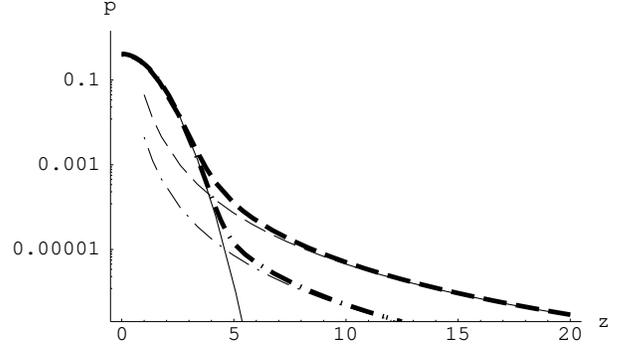}} \vspace*{10mm} \caption{\protect{Plots of
Gaussian distribution Eq.~(\ref{GD}) shown as a thin full line and
probability densities $p_{_{N}}(z)$ calculated numerically
according to Eqs.~(\ref{pl>2}), (\ref{gsum}) at $l=2$ for $N=10^2$
(thick dashed) and $N=10^4$ (thick dot-dashed). Both $p_{_N}(z)$
exhibit tails {\it heavier\/} than the Gaussian. The corresponding
tail asymptotics, Eq.~(\ref{plf}) are drown as thin dashed and
dot-dashed lines respectively. }\label{F3}}
\end{figure}

Comparing Eq.~(\ref{plf}) with Eq.~(\ref{p1as}) one can see that
the two expressions yield exactly the same asymptotic behavior for
the probability densities\cite{note}, including the values of the
normalizing multiplier [see the definition of $A$ in
Eq.~(\ref{fl})]. It should be stressed, however, that this is not
the case for any $f(\xi)$ with a tail lighter the the Gaussian.
Eq.~(\ref{p1as}) describes the probability density for a sum to
have a great value because of contribution of a single extremely
great term. Meanwhile if $f(\xi)$ decays sharp enough the
probability of such an event is very small. In this case it could
be more probable to obtain the great value of the sum due to
contribution of {\it several\/} great but not extremely great
terms. It means, that the actual tail of $p_{_N}(x)$ just cannot
be {\it lighter\/} than that given by Eq.~(\ref{p1as}) but, in
principle, it could be {\it heavier}. To elucidate this issue let
us consider the following profile of $f(\xi)$
\[  f(\xi) = \frac{1}{\pi \cosh \xi}.
\]
The corresponding variance $\sigma^2$ equals $\pi^2/4$. Function
$g_{_N}(\omega)$ is calculated according to Eq.~(\ref{gN}) with
standard methods of integration in a complex plane. The
integration yields
\[  g_{_N}(\omega) = \frac{1}{\cosh{\frac{\omega \pi}{2\sigma \sqrt
  N}}} \equiv \frac{1}{\cosh{\frac{\omega}{\sqrt N}}}.
\]
Then, $p_{_N}(z)$ is given by the following integral
\begin{equation}\label{pch}
  p_{_N}(z) = \frac{\sigma\sqrt N}{\pi^2}\int_{-\infty}^{\infty}
  \frac{\exp(-2iz\zeta\frac{\sigma\sqrt N}{\pi})}{\cosh^N \zeta}d\zeta,
\end{equation}
where $\zeta \equiv \omega \pi/(2\sigma \sqrt N)$. To find the
asymptotics of integral Eq.~(\ref{pch}) at $z^2 \gg \ N$ let us
extend the integration to the plane of complex $\zeta$ and
consider a standard contour, where the edges of the segment $-R
\leq \zeta^\prime \leq R$ are connected through the lower
half-plane by a circular arc with the radius $R$. The integrand
has infinite number of poles of order $N$ lying in points
$\zeta^{\prime\prime} = (2n+1)\pi/2$ of the imaginary axis. Here
$n$ is any integer. Bearing in mind that the integral over the arc
tends to zero at $R \rightarrow \infty$ and the direction of
circulation about the contour is negative one obtains
\begin{eqnarray}
 & & p_{_N}(z) =  \label{res} \\
 & & -2\pi i \frac{\sigma\sqrt N}{\pi^2}
  \sum_{n=0}^{\infty}{\rm res}\left[ \frac{\exp(-2iz\zeta
  \frac{\sigma\sqrt N}{\pi})}{\cosh^N \zeta};\;
  -i\pi\frac{2n+1}{2}\right].\nonumber
\end{eqnarray}
In its turn
\begin{eqnarray}
     &  & {\rm res}\left[ \frac{\exp(-2iz\zeta
  \frac{\sigma\sqrt N}{\pi})}{\cosh^N \zeta};\;
  -i\pi\frac{2n+1}{2}\right] = \nonumber \\
     &  & \hspace*{17mm}
     \lim_{\zeta \rightarrow -i\pi(2n+1)/2}
     \left\{\frac{1}{(N-1)!}
     \frac{d^{N-1}}{d\zeta^{N-1}}\label{res1}\right.\\
     &  &\left. \left[\left(
     \frac{\zeta + i\pi\left(n+\frac{1}{2}\right)}
     {\cosh \zeta}\right)^N \!\!\!
     \exp\left(-2iz\zeta
     \frac{\sigma\sqrt N}{\pi}\right)
     \right]\right\}. \nonumber
  \end{eqnarray}

We are interested in the leading (in $z \gg \sqrt N$) term of
Eq.~(\ref{res1}). Note, that each differentiation of $(\ldots)^N$
contributes a multiplier of order of $N$, while each
differentiation applied to $\exp(\ldots)$ yields the one of order
of $z\sqrt N \gg N$. Thus, to obtain the leading term of
Eq.~(\ref{res1}) one must apply all $N-1$ derivatives to the
exponent. It gives rise to the following expression for the
residue
\begin{eqnarray}
&  & {\rm res}\left[ \frac{\exp(-2iz\zeta
  \frac{\sigma\sqrt N}{\pi})}{\cosh^N \zeta};\;
  -i\pi\frac{2n+1}{2}\right] \cong \label{resf} \\
& & i\frac{(-1)^{nN}}{(N-1)!}\left(2z\frac{\sigma \sqrt
    N}{\pi}\right)^{N-1}\!\!\!\exp\left[-(2n+1)z\sigma \sqrt N\right]
    +..,\nonumber
\end{eqnarray}
where dots denote dropped higher order (in $\sqrt N/z$) terms.

Substitution of Eq.~(\ref{resf}) in Eq.~(\ref{res}) transforms the
sum there into a geometric progression, which may be summarized.
However, the summarizing is not required --- due to the
exponential decay of the sum's terms the sum equals to its first
term with the exponential accuracy. Finally the desired tail of
the probability density at $z \gg \sqrt N$ is given by the
following formula
\begin{eqnarray}\label{pchf}
p_{_N}(z) \cong \frac{(2\sigma\sqrt N/\pi)^N}{(N-1)!} & z^{N-1} &
\exp\left(-z\sigma \sqrt
  N\right) \equiv \nonumber\\
   \frac{N^{N/2}}{(N-1)!} & z^{N-1} & \exp\left(-\frac{\pi z \sqrt
  N}{2}\right),
\end{eqnarray}
(I remind to the reader that $\sigma = \pi/2$). It should be
emphasized that the above expression gives rise to $p_{_N}(z)$,
which is much greater than that following from Eq.~(\ref{p1as}). A
plot of $p_{_N}(z)$ obtained by numerical integration of
Eq.~(\ref{pch}) as well as its comparison with GD and asymptotic
Eq.~(\ref{pchf}) are shown in Fig.~\ref{F4}.
\begin{figure}
\epsfxsize = 0.45\textwidth 
\makebox{
\epsfbox{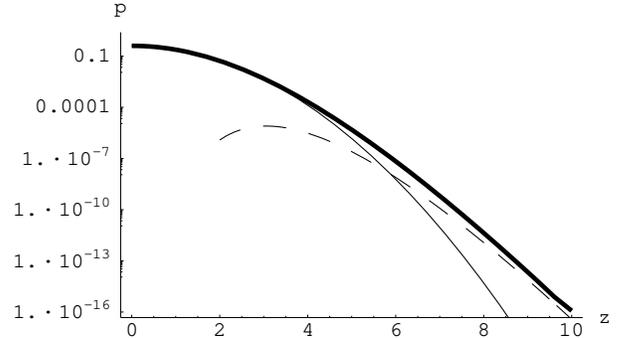}} \vspace*{10mm} \caption{\protect{Plots of
probability density $p_{_{N}}(z)$ calculated numerically according
to Eq.~(\ref{pch}) at $N=25$ (thick line), Gaussian distribution
Eq.~(\ref{GD}) for the normalized sum $z$ (thin line) and
asymptotics Eq.~(\ref{pchf}) shown as a dashed line. See the text
for details.}\label{F4}}
\end{figure}

Thus, the above consideration has shown that the very tail of the
probability density $p_{_N}(x)$ for a sum $x$ of $N$ random terms
$\xi_n$ at any finite $N$ is defined by the tail of the
probability density $f(\xi)$. The tail is lighter than a Gaussian,
if $f(\xi)$ decays at large $\xi$ faster than GD, and heavier than
it in the opposite case. The only situation, when the tail of
$p_{_N}(x)$ is the Gaussian, corresponds to a Gaussian $f(\xi)$.
For $f(\xi)$ decaying as a certain power of $\xi$ large
fluctuations of $x$ occur because of contribution to the sum of a
single anomalously great term, rather than due to the ones of a
number of terms of moderate values\cite{note1}, while for
exponentially decaying $f(\xi)$ this is not the case.

I benefited much from discussions of this work with N. Makarenko
and H. Takayasu. This study was supported by the Grant-in-Aid for
Scientific Research (No. 11837006) from the Ministry of Education,
Culture, Sports, Science and Technology (Japan).

\appendix
\section*{}
Here I discuss the derivation of Eq.~(\ref{plf}). Note, that if
$\omega$ is regarded as a complex variable, then $g_{_N}(\omega)$
given by Eq.~(\ref{gsum}) is an analytical function of it at any
finite $|\omega|$. Since $N$ is a positive integer number the same
is true for $\left(g_{_N}(\omega)\right)^N$. For this reason
$\int_C \left(g_{_N}\right)^N\exp(-i\omega z)d\omega = 0$, where
$C$ stands for any closed contour in the plane of complex
$\omega$.

To obtain the asymptotic behavior of $p(z)$ at large $z$ let us
consider contour $C$ consisting of segments of real ($0 \leq
\omega^\prime \leq R$ and imaginary ($0 \geq \omega^{\prime\prime}
\geq -R$) axes and a circular arc with the radius $R$, connecting
the edges of the segments. According to Eq.~(\ref{gsum}) the
integrand decays along the real axis as $\exp\left[-(\omega^\prime
\sqrt N/\sigma)\sin(\pi/2l)\right]$. In the lower half-plane it
decays not slower than $\exp\left[\left(|\omega|\sqrt
N/\sigma\right)-\omega^{\prime\prime}z\right]$. Then, it is easy
to show that at $z \gg \sqrt N/\sigma$ the integral over the
circular arc tends to zero at $R\rightarrow \infty$. Thus, in the
specified limit $z \gg \sqrt N/\sigma$ one obtains
\begin{eqnarray}
   & p_{_N}(z)& \cong \nonumber \\
   & & -\frac{1}{\pi}{\rm Re}\int_{-i\infty}^0
  \left(g_{_N}(\omega)\right)^N
  e^{-i\omega z}d\omega
  \equiv \label{pNz}\\
 & & \frac{1}{\pi}{\rm Im}\int_0^\infty
  \left(g_{_N}(-i\omega^{\prime\prime})\right)^N
  e^{-\omega^{\prime\prime}z}d\omega^{\prime\prime}\nonumber
\end{eqnarray}
Next, note that $\exp(-\omega^{\prime\prime}z)$ is a sharp
function relative to
$\left(g_{_N}(-i\omega^{\prime\prime})\right)^N$. In this case to
find the leading asymptotic of $p_{_N}(z)$ at $z \rightarrow
\infty$ one may expand $g_{_N}(-i\omega^{\prime\prime})$ in powers
of $\omega/(\sigma \sqrt N)$, truncating the expansion when the
first term with nontrivial contribution to $p_{_N}(z)$ is
obtained.

The replacement in Eq.~(\ref{gsum}) $\omega$ with
$-i\omega^{\prime\prime}$ and the expansion of exponents in the
numerators results in 
\begin{eqnarray}
 & &   g_{_N}(-i\omega^{\prime\prime}) =\nonumber\\
     &  &i\left( \sin \frac{\pi}{2l} \right)\sum_{j=o}^{l-1}\sum_{k=0}^\infty
     \frac{1}{k!}\left(\frac{\omega^{\prime\prime}}{\sigma \sqrt N}\right)^k
     \times     \nonumber\\
     & & \hspace*{2cm}
     \exp\left[\frac{i\pi}{2l}(2j+1)(k-2l+1)\right]\label{fl>2}\\
    & & \equiv -i\left( \sin \frac{\pi}{2l} \right)\sum_{k=0}^\infty
    \frac{1}{k!}\left(\frac{\omega^{\prime\prime}}{\sigma \sqrt N}\right)^k
    \exp\left[ \frac{i\pi}{2l}(k+1)\right]\times\nonumber\\
    &  & \hspace*{2cm}
    \sum_{j=0}^{l-1}\exp\left[\frac{i\pi(k+1)}{l}j\right]\nonumber
  \end{eqnarray}
The latter sum over $j$ is a geometric progression. Summarizing it
up one obtains
\begin{eqnarray*}
 \sum_{j=0}^{l-1} & \exp & \left[\frac{i\pi(k+1)}{l}j\right]=\\
& & \left\{
  \begin{array}{l}
    2/\left[1 - \exp\left(i\pi \frac{k+1}{l}\right)\right] \;  {\rm
    at}\; k=2p\\
    0\;  {\rm at}\; k=2p+1 \neq 2ql - 1 \\
    l\;  {\rm at}\; k=2ql-1
  \end{array}
\right.
\end{eqnarray*}
Here $p=0,\;1,\;2,\ldots;\; q=1,\;2,\;3,\ldots$

Then,  Eq.~(\ref{fl>2}) may be transformed as follows
\begin{equation}\label{rs}
 g_{_N}(-i\omega^{\prime\prime}) =
 \sum_{p=0}^{\infty} r_p +  i\sum_{q=1}^\infty s_q \;,
\end{equation}
where
\begin{eqnarray*}
     &  & r_p =
     \frac{1}{(2p)!}
     \left(\frac{\omega^{\prime\prime}}
     {\sigma \sqrt N}\right)^{2p}
     \frac{\sin \frac{\pi}{2l}}{\sin
     \frac{\pi(2p+1)}{2l}} \;\;,\\
     &  & s_q =
     \left(l \sin\frac{\pi}{2l}\right)
     \frac{(-1)^{q-1}}{(2ql-1)!}
     \left(\frac{\omega^{\prime\prime}}
     {\sigma \sqrt N}\right)^{2ql-1}.
  \end{eqnarray*}

At small $\omega$, retaining just the first terms in both the sums
in Eq.~(\ref{rs}) one arrives at the following result
$g_{_N}(-i\omega^{\prime\prime}) \cong 1 + ... + i s_1 + ..,$ and
next, $ {\rm Im}\left( g_{_N}(-i\omega^{\prime\prime})\right)^N
\cong Ns_1 + ..,$ where dots denote dropped higher order in
$\omega^{\prime \prime}/(\sigma \sqrt N)$ terms. Substitution it
into Eq.~(\ref{pNz}) and integration yield Eq.~(\ref{plf}).

\end{twocolumns}

\begin{references}
\vspace*{-10mm}
\bibitem[*]{E-mail}E-mail: tribel@scroll.apphy.fukui-u.ac.jp
\bibitem{GK} B. V. Gnedenko, and A. Kolmogorov, {\it Limit
Distribution for Sums of Independent Random Variables}
(Addison-Wesley, Cambridge, MA, 1954).
\bibitem{SF}G. Boffetta, V. Carbone, P. Giuliani {\it et al},
Phys. Rev. Lett. {\bf 83}, 4662(1999).
\bibitem{war}D. Roberts, and D. L. Turcotte Fractals, {\bf 6}, 351
(1998).
\bibitem{Stan}V. Plerou, P. Gopikrishnan, L. A. Nunes Amaral {\it et
al}, Phys. Rev. E {\bf 60}, 6519 (1999).
\bibitem{turb}T. Bohr, M. H. Jensen, G. Paladin and A. Vulpiani,
{\it Dynamical Systems Approach to Turbulence} (Cambridge
University Press, Cambridge U.K. 1998).
\bibitem{Tak}Taking the opportunity, I would like to thank H.
Takayasu, who kindly pointed me out this simple and beautiful
example of a non-Gaussian tail.
\bibitem{note}For the correct comparison of the probability
densities one should remember that change of a random variable
results in the change of the corresponding probability density
according to the formula $p(z)=p(x(z))(dx/dz)$.
\bibitem{n}Note, to provide statistical description for $x$ one
should consider an ensemble of sequences $\{\xi_n\},\; 1 \leq n
\leq N$ with the same value of $N$ for any sequence in the
ensemble. This approach does not depend on the particular $N$.
Thus, not only large but also {\it any\/} positive integer values
of $N$ are meaningful.
\bibitem{Zas}See, e.g., M. F. Shlesinger, G. M. Zaslavsky, and J.
Klafter Nature {\bf 363}, 31 (1993).
\bibitem{note1}Strictly speaking the latter was proved just for
$f(\xi)$ decaying as even power $\xi$. However, bearing in mind
Eq.~(\ref{m}) the extension of the conclusion to any power decay
sounds quite reasonable.
\end{references}
\end{document}